\documentclass[12pt,twoside]{article}
\usepackage{amsfonts}
\usepackage{amssymb,amsmath}
\usepackage{setspace}
\begin{document}
\begin{center} \textbf{Note on the Radicals of Numbers} \vspace{12pt}
\\Constantin M. Petridi
\\ cpetridi@hotmail.com\\
\vspace{20pt}
 \textbf{Abstract}\end{center}
\par
\hspace{0.2in}
\begin{tabular}{p{12cm}}
\begin{small}We prove that for any given $\theta > 0$ the
geometric mean of the radicals of the $\phi(n)$ integers less
than $n$ and relatively prime to $n$, is greater than
$K(\theta)\;R(n)^{-\theta}\;n$, where $K(\theta)$ an absolute
constant depending only on $\theta$.
\end{small}
\end{tabular}
\vspace{25pt}
\par
Denoting by $R(n)$ the radical $p_{1}p_{2}\ldots p_{\omega}$ of
any positive integer $n = p_{1}^{a_{1}}p_{2}^{a_{2}}\ldots
p_{\omega}^{a_{\omega}}$, let $$R(n_{1}),R(n_{2}),\ldots
,R(n_{\phi(n))}$$ be the radicals of the $\phi(n)$ integers less
than $n$ and relatively prime to $n$.
\par
We prove following theorem for the geometric mean
$$\left ( \prod_{i=1}^{\phi(n)}R(n_{i}) \right
)^{\frac{1}{\phi(n)}}.$$
\par
 \textbf{Theorem.} Given any $\theta >0$, $$\left ( \prod_{i=1}^{\phi(n)}R(n_{i}) \right
)^{\frac{1}{\phi(n)}} > K_{\theta}\;R(n)^{-\theta}\;n,$$ where
$K_{\theta}$ an absolute constant depending only on $\theta$.
\par
\textbf{Proof.} In our paper \textit{A strong "abc conjecture"
for certain partitions $a+b$ of $c$} (arXiv..math.NT/0511224 v2 5
Jan 2006) we proved (Theorem 3) that for any given $\varepsilon
>0$
$$\hspace{4cm}G_{n}^{\frac{2}{\phi(n)}}>k_{\varepsilon}\;R(n)^{1-\varepsilon}\;n^{2},
\hspace{5cm}(1)$$ where
$$G_{n} =\prod_{i=1}^{\frac{\phi(n)}{2}}R(x_{i}y_{i}n)$$
extends over all $\phi(n)/2$ positive solutions $x_{i}, y_{i}$ of
$$x_{i}+y_{i}=n \hspace{30pt} i = 1 \ldots \phi(n)/{2},$$ with
$x_{i}<y_{i}$ and $(x_{i},y_{i})=1$.
\par
Dividing both sides of (1) by $R(n)$, and considering that
$R(x_{i}y_{i}n)=R(x_{i})R(y_{i})R(n)$, we obtain $$\left (
\prod_{i=1}^{\frac{\phi(n)}{2}}R(x_{i})R(y_{i}) \right
)^{\frac{2}{\phi(n)}}
> k_{\varepsilon}\;R(n)^{-\varepsilon}\;n^{2}.$$
which raised to the power $1/2$, gives
$$\hspace{3.5cm}\left (
\prod_{i=1}^{\frac{\phi(n)}{2}}R(x_{i})R(y_{i}) \right
)^{\frac{1}{\phi(n)}}
>\; k_{\varepsilon}^{\frac{1}{2}}\;R(n)^{-\frac{\varepsilon}{2}}\;n. \hspace{4cm} (2)$$
As the set of the $\phi(n)$ numbers $x_{i},y_{i}, \hspace{5pt} i
= 1, \ldots \phi(n)/2$, is identical with the set $n_{1},n_{2},
\ldots n_{\phi(n)}$, defined above, (2) can be written as
$$\left (
\prod_{i=1}^{\phi(n)}R(n_{i} \right )^{\frac{1}{\phi(n)}}
> \;k_{\varepsilon}\;R(n)^{-\frac{\varepsilon}{2}}\;n.$$
Setting $\varepsilon/2=\theta>0$ and
$k_{\varepsilon}^{1/2}=k_{\theta}$ as a new absolute constant
depending only on $\theta$ we obtain
$$\left (
\prod_{i=1}^{\phi(n)}R(n_{i} \right )^{\frac{1}{\phi(n)}}
> \;k_{\theta}\;R(n)^{-\theta}\;n,$$ as claimed by the theorem.
\end{document}